\ifx\LocalElevenPtArticleFormatLoaded\undefined
  \documentclass[a4paper,11pt]{article}     
  \usepackage{array}                        
  \usepackage{theorem}                      
  \usepackage{amsmath,amscd,amssymb}        
  \usepackage{latexsym}                     
  \usepackage[german,english]{babel}        
  \usepackage[applemac]{inputenc}
 \usepackage{xypic}
\fi

\newtheorem{theorem}{Theorem}[section]
\usepackage{theorem}

\newtheorem{corollary}[theorem]{Corollary}

\theorembodyfont{\rm}
\newtheorem{remark}{Remark}

\newcommand{\CC}{{\mathbb{C}}}

\newcommand{\HH}{{\mathbb{H}}}

\newcommand{\PP}{{\mathbb{P}}}
\newcommand{\QQ}{{\mathbb{Q}}}

\newcommand{\ZZ}{{\mathbb{Z}}}

\newcommand{\op}{\operatorname}

\newenvironment{Proof}{\begin{ProofwCaption}{Proof}}{\end{ProofwCaption}}
\newenvironment{Proof*}[1]{\begin{ProofwCaption}{{#1}}}{\end{ProofwCaption}}
\newenvironment{ProofwCaption}[1]%
  {\addvspace\theorempreskipamount \noindent{\it #1.}\rm}%
  {\qed \par \addvspace\theorempostskipamount}
\newcommand{\qedsymbol}{\mbox{$\Box$}}
\newcommand{\qed}{\quad\qedsymbol}
\begin{document}

\title{A remark on the Schottky locus in genus $4$}
\author{ Joe Harris and Klaus Hulek}
\date{}
\maketitle

\section{Introduction}

By Torelli's theorem the map $t_g:{\cal M}_g\rightarrow {\cal A}_g$ which
associates to each curve $C$ its Jacobian variety $J(C)$, together with the
principal polarization given by the theta divisor, is injective. We denote its
image by $J_g$. In genus $4$ Mumford \cite{Mu} has determined the class of the
closure $J'_4$ of $J_4$ in the partial compactification ${\cal A}'_4$ which is
given by adding the rank 1 degenerations. In this brief note we shall determine
the class of the Schottky locus in the Igusa and the Voronoi compactifications
of ${\cal A}_4$ and we shall briefly comment on the degree $8$ 
modular form which vanishes on the Jacobi locus $J_4$.

\section{The main result} 
We denote by ${\cal A}_g^{\op{Igu}}$ and ${\cal A}_g^{\op{Vor}}$ the Igusa,
resp. the second Voronoi compac\-ti\-fi\-cation of ${\cal A}_g$. 
These are toroidal compactificiations given by choosing either the
second Voronoi decomposition (which gives the second Voronoi compactification) 
or the central cone decomposition, which
for genus $g=4$ coincides with the perfect cone decomposition (and which
leads to the Igusa compactification).
If $g=4$, then
${\cal A}_4^{\op{Igu}}$ is singular and ${\cal A}_4^{\op{Vor}}$ is smooth in
the stack sense. By this we mean the following: the varieties
${\cal A}_4^{\op{Vor}}$ and ${\cal A}_4^{\op{Igu}}$  are both singular, but
${\cal A}_4^{\op{Vor}}$ only has finite quotient singularities and adding a
full
level-$n$ structure of level $n\ge 3$ makes ${\cal A}_4^{\op{Vor}}(n)$
smooth, whereas ${\cal A}_4^{\op{Igu}}(n)$ will always be singular. 
These singularities come from the choice of the fan.
As a stack
${\cal A}_4^{\op{Igu}}$ has exactly  one singularity (and by abuse of 
notation we will refer to this as the singularity of  ${\cal A}_4^{\op{Igu}}$. 
There is a morphism
$\pi:{\cal A}_4^{\op{Vor}}\rightarrow {\cal A}_4^{\op{Igu}}$ which is a blow-up
with centre at the singular point of ${\cal A}_4^{\op{Igu}}$ (see
\cite[chapter I]{HS}). (Note that this is special for genus $4$.) 
We denote by $D'_4$, resp. $D_4^{\op{Igu}}$ and
$D_4^{\op{Vor}}$ the boundary component given by the rank 1 degenerations,
resp. its closure in ${\cal A}_4^{\op{Igu}}$ and ${\cal A}_4^{\op{Vor}}$ and we
denote by $E$ the (reduced) irreducible divisor in ${\cal A}_4^{\op{Vor}}$
which
is exceptional with respect to
$\pi:{\cal A}_4^{\op{Vor}}\rightarrow {\cal A}_4^{\op{Igu}}$.
Finally, let $L$
be the $(\QQ)$-line bundle of modular forms of weight 1 on the Satake
compactification ${\cal A}_4^{\op{Sat}}$. Since both ${\cal A}_4^{\op{Igu}}$
and ${\cal A}_4^{\op{Vor}}$ map to the Satake compactification we can also
consider $L$ as a line bundle on these varieties, which, by abuse of notation,
we will again denote by $L$. 
We denote by ${\cal A}'_4$ Mumford's partial compactification which is
given by adding the rank $1$ degenerations of principally polarized abelian
varieties. Note that ${\cal A}_4^{\op{Igu}} \setminus {\cal A}'_4$ has 
codimension $2$ whereas  ${\cal A}_4^{\op{Vor}} \setminus {\cal A}'_4$ contains
the exceptional divisor $E$.
We finally recall from \cite[Prop. I.1 and
Prop. I.6]{HS} that the inclusion $i:{\cal A}'_4\rightarrow {\cal
A}_4^{\op{Igu}}$ induces an isomorphism
$$
i^{\ast}:\mbox{Pic } {\cal A}_4^{\op{Igu}}\otimes\QQ=\QQ L\oplus\QQ
D_4^{\op{Igu}}\cong \mbox{Pic } {\cal A}'_4 \otimes\QQ
$$
and that
$$
\mbox{Pic } {\cal A}_4^{\op{Vor}}\otimes \QQ=\QQ L\oplus\QQ
D_4^{\op{Vor}}\oplus
\QQ E.
$$
Moreover by \cite[Prop I.3]{HS}
$$
\pi^{\ast} D_4^{\op{Igu}}=D_4^{\op{Vor}}+4E.
$$
The Torelli map $t_g:{\cal M}_g\rightarrow {\cal A}_g$ extends by \cite{N} to a
morphism $\bar{t}_g:{\overline{\cal M}}_g\rightarrow {\cal A}_g^{\op{Vor}}$ on
the Deligne-Mumford compactification of ${\cal M}_g$. We denote its image by
$J_4^{\op{Vor}}$ and we set $J_4^{\op{Igu}}=\pi \left(J_4^{\op{Vor}}\right)$.
This coincides with the closure of $J'_4$ in ${\cal A}_4^{\op{Igu}}$.

\begin{theorem} \label{theorem}
The varieties ${\cal A}_4^{\op{Vor}}$ and ${\cal A}_4^{\op{Igu}}$ are
$\mathbb Q$-factorial varieties. For the Schottky locus the 
following holds:
\begin{enumerate}
\item[\rm{(i)}]
The Schottky locus $J_4^{\op{Vor}}$ does not meet the exceptional divisor $E$
in ${\cal A}_4^{\op{Vor}}$ and hence $J_4^{\op{Igu}}$ does not pass through the
singular point of ${\cal A}_4^{\op{Igu}}$.
\item[\rm{(ii)}]
The divisor $J_4^{\op{Igu}}$ is $\mathbb Q$-Cartier in  
${\cal A}_4^{\op{Igu}}$ and its class equals
$$
[J_4^{\op{Igu}}] = 8L - D_4^{\op{Igu}}.
$$
\item[\rm{(iii)}]
The divisor $J_4^{\op{Vor}}$ is $\mathbb Q$-Cartier in ${\cal A}_4^{\op{Vor}}$ 
and its class equals
$$
[J_4^{\op{Vor}}] = 8L - D_4^{\op{Vor}}-4E.
$$
\end{enumerate}
\end{theorem}

\begin{Proof}
We shall first prove that  ${\cal A}_4^{\op{Vor}}$ 
and ${\cal A}_4^{\op{Igu}}$ are
$\mathbb Q$-factorial varieties. This is clear for
${\cal A}_4^{\op{Vor}}$ since it is the quotient of the smooth
variety ${\cal A}_4^{\op{Vor}}(n)$, for $n \geq 3$ by a finite group,
namely ${\op{PSL(8,n)}}$. In particular also ${\cal A}'_4$ is 
$\mathbb Q$-factorial and by Mumford's result we know that 
$ \mbox{Pic } {\cal A}'_4 \otimes\QQ$ is generated by $L$ and
$D'_4$. Both these divisors extend to $\mathbb Q$-Cartier divisors
on ${\cal A}_4^{\op{Igu}}$ where they generate the group 
$\mbox{Pic } {\cal A}_4^{\op{Igu}}\otimes\QQ$.

Next we shall prove that $J_4^{\op{Vor}}$ and the exceptional locus $E$ have
empty intersection. Assume that this is not the case. Since $E$ is a
$\QQ$-Cartier divisor on ${\cal A}_4^{\op{Vor}}$ this would imply that $(\bar
t_4)^{-1}(E)$ has codimension $1$. Recall that the points on $E$ all correspond
to degenerations of abelian varieties with trivial abelian part. This would
imply (cf.\cite{A})  that $(\bar t_4)^{-1}(E)$ is a 
$\mathbb Q$-divisor whose
general point 
for each component
corresponds to a stable curve whose normalization has only
rational components. But this is impossible since ${\cal M}_4$ maps to ${\cal
A}_4$ and the boundary
$\overline{\cal M}_4\setminus {\cal M}_4$ consists of 3 components ${\cal
M}_{4-i, i}; i=0,1,2$ whose general points consist of a curve of genus $3$ with
one node, a curve of genus 3 and an elliptic curve intersecting in one point,
or two genus 2 curves intersecting in one point respectively. These are
$8$-dimensional families, but their general point is mapped to $D'_4$ in the
first case, resp. the interior ${\cal A}_4$ in the other two cases.
In particular, the Jacobian of the general point of such a divisor has an
abelian part of dimension either $3$ or $4$.

It was shown by Mumford \cite[p.369]{Mu} that
$$
[J'_4]=8L-D'_4
$$
and this, together with the fact that the inclusion $i:{\cal A}'_4\rightarrow
{\cal A}_4^{\op{Igu}}$ induces an isomorphism of the Picard 
groups, proves (ii).
Since $\pi^{\ast} D_4^{\op{Igu}}=D_4^{\op{Vor}}+4E$ the class of the total
transform of $J_4^{\op{Igu}}$ in ${\cal A}_4^{\op{Vor}}$ equals
$\pi^{\ast}(J_4^{\operatorname{Igu}})=
8L-D_4^{\op{Vor}}-4E$. But since $J_4^{\op{Igu}}$ does not pass through the
singular point of ${\cal A}_4^{\op{Igu}}$ this also proves (iii).
\hfill\end{Proof}

\section{Modular forms}
There exists a weight 8 modular form $F$ which vanishes precisely along $J_4$.
There are two ways to write this down explicitly. The first goes back to Igusa.
Let $S$ be an even, integer, positive definite quadratic form with
determinant 1 of rank $m$. Such a form defines a genus $g$ modular form of
weight
$m/2$ given by the theta series
$$
\vartheta^{(g)}_s(\tau)=\sum\limits_{G\in \operatorname{ Mat}(m,g;\ZZ)}\exp[\pi
i\mbox{Tr}(^tGSG\tau)].
$$
There is only one even, integer, positive definite quadratic form $S(8)$ with
determinant 1 of rank 8 and there are two such forms $S(8)\oplus S(8)$ and
$S(16)$ of rank 16. Igusa \cite{I}, cf. also \cite[p.20]{W} has shown that
$$
F(\tau)=\vartheta^{(4)}_{S(8)\oplus S(8)} (\tau)-\vartheta^{(4)}_{S(16)}(\tau)
$$
is a non-zero cusp form of weight 8 which vanishes on $J_4$.\\
The second method goes back to van Geemen and van der Geer. 
For $\varepsilon, \varepsilon'\in\{0,1\}^{g-1}$  define the following theta
functions
$$
\theta [\begin{matrix}\varepsilon\\[-2mm] \varepsilon'\end{matrix}]
(\tau,z)=\sum_{m\in\ZZ^{g-1}} \mbox{exp}
[\pi i{}^t((m+\frac{\varepsilon}{2})\tau(m+\frac{\varepsilon}{2})+2{}^t
(m+\frac{\varepsilon}{2})(z+\frac{\varepsilon'}{2}))],
$$
resp. the theta constants
$$
\theta[\begin{matrix}\varepsilon\\[-2mm] \varepsilon'\end{matrix}]=
\theta[\begin{matrix}\varepsilon\\[-2mm] \varepsilon'\end{matrix}](\tau,0).
$$
Such a theta function is called {\em even} (resp. {\em odd}) if
$^t\varepsilon\varepsilon'\equiv 0\mod 2$ (resp.
$^t\varepsilon\varepsilon'\equiv 1\mod 2$). The odd theta constants vanish
identically and there are $2^{g-2}(2^{g-1}+1)$ even theta constants. They
define a morphism
$$
\begin{array}{rcl}
\operatorname {Sq} : {\cal A}_{g-1}(2,4) &\longrightarrow & \PP^N, \quad
N=2^{g-2}(2^{g-1}+1)-1\\
{[}\tau{]}  &\longmapsto & (\ldots,
\theta^2[\begin{matrix}\varepsilon\\[-2mm] \varepsilon'\end{matrix}](\tau,0),
\ldots).
\end{array}
$$
Van Geemen \cite{vG1}, resp. van Geemen and van der Geer \cite{vGvG} have shown
how elements in the ideal of $\operatorname{Sq}({\cal A}_{g-1}(2,4))$ give
rise to
modular forms vanishing on the Schottky locus  $J_g$ in ${\cal A}_g$. This
can be
described as follows. Let $\CC[X_0,\ldots,X_N]$ be the coordinate ring of
$\PP^N$, and let $I_{g-1}$ be the homogeneous ideal of
$\operatorname{Sq}({\cal A}_{g-1}(2,4))$.
Then for every element $G(X_0,\ldots, X_N)$ of $I_{g-1}$ we define a
form on
$\HH_g$ as follows
$$
\sigma(G)(\tau)=G(\ldots,\theta\left[\begin{matrix}\varepsilon & 0\\
\varepsilon' &
0\end{matrix}\right](\tau,
0)\theta\left[\begin{matrix}\varepsilon & 0\\ \varepsilon' &
1\end{matrix}\right] (\tau, 0),\ldots); \quad \tau \in \HH_g.
$$
This will then be a modular form vanishing on $J_g$.

In the case of genus $g=4$ we, therefore, have to start with a theta relation
for $g=3$. According to \cite[p. 623]{vGvG} one has the following relation
$$
\begin{array}{l}
{\ \ \, }\theta\left[\begin{matrix}0&0&0\\0&0&0\end{matrix}\right]
\theta\left[\begin{matrix}0&0&0\\1&0&0\end{matrix}\right]
\theta\left[\begin{matrix}0&0&0\\0&1&0\end{matrix}\right]
\theta\left[\begin{matrix}0&0&0\\1&1&0\end{matrix}\right]\\[4mm]
-\theta\left[\begin{matrix}0&0&1\\0&0&0\end{matrix}\right]
\theta\left[\begin{matrix}0&0&1\\1&0&0\end{matrix}\right]
\theta\left[\begin{matrix}0&0&1\\0&1&0\end{matrix}\right]
\theta\left[\begin{matrix}0&0&1\\1&1&0\end{matrix}\right]\\[4mm]
-\theta\left[\begin{matrix}0&0&0\\0&0&1\end{matrix}\right]
\theta\left[\begin{matrix}0&0&0\\1&0&1\end{matrix}\right]
\theta\left[\begin{matrix}0&0&0\\0&1&1\end{matrix}\right]
\theta\left[\begin{matrix}0&0&0\\1&1&1\end{matrix}\right]=0
\end{array}
$$
If one writes this relation as $r_1-r_2-r_3=0$, one obtains from this the
following relation between the squares
$\theta^2\left[\begin{matrix}\varepsilon\\[-2mm]
\varepsilon'\end{matrix}\right]$:
$$
r^4_1+r^4_2+r^4_3-2r^2_1r^2_2-2r^2_1r^2_3-2r^2_2r^2_3=0.
$$
Replacing $\theta^2\left[\begin{matrix}\varepsilon\\[-2mm]
\varepsilon'\end{matrix}\right]$ by
$\theta\left[\begin{matrix}\varepsilon&0\\[-2mm]
\varepsilon'&0\end{matrix}\right]
\theta\left[\begin{matrix}\varepsilon&0\\[-2mm]
\varepsilon'&1\end{matrix}\right]$ then yields a modular form on $\HH_4$ of
weight $8$. Up to a possibly non-zero constant this is equal to $F$.
We can consider $F$ as a section of the $\mathbb Q$-line bundle $L$
both on  ${\cal A}_4^{\op{Igu}}$ and ${\cal A}_4^{\op{Vor}}$. We conclude 
this note with a remark on the divisor defined by $F$, which in the case of the
Igusa compactification is essentially nothing but Mumford's observation.

\begin{corollary}
\begin{enumerate}
\item[\rm{(i)}]
The divisor of $F$ on ${\cal A}_4^{\op{Igu}}$ equals
$J_4^{\op{Igu}}+D_4^{\op{Igu}}$.
\item[\rm{(ii)}]
The divisor of $F$ on ${\cal A}_4^{\op{Vor}}$ equals
$J_4^{\op{Vor}}+D_4^{\op{Vor}}+4E$.
\end{enumerate}
\end{corollary}

\begin{Proof}
The divisor defined by $F$ on ${\cal A}'_4$ equals $J'_4+D'_4\sim 8L$. (This
gives Mumford's claim that $J'_4=8L-D'_4.$) Since
${\cal A}'_4$ and ${\cal A}_4^{\op{Igu}}$ only differ in codimension 2 the
analogous statement is true on ${\cal A}_4^{\op{Igu}}$. Moreover, since
$\pi^{\ast} D_4^{\op{Igu}}=D_4^{\op{Vor}}+4E$ the form $F$, considered as a
form
on ${\cal A}_4^{\op{Vor}}$ vanishes at least on
$J_4^{\op{Vor}}+D_4^{\op{Vor}}+4E$. On the other hand we know that the class of
$J_4^{\op{Vor}}$ equals $8L-D_4^{\op{Vor}}-4E$ by  theorem \ref{theorem}
and hence $F$
vanishes on $E$ of order exactly $4$.
\hfill\end{Proof}

\begin{remark}
It should be possible to compute the vanishing order
of $F$ on $E$ from the Fourier expansion of
$F$ with respect to $E$, but this seems to be a cumbersome computation.
\end{remark}

\bibliographystyle{alpha}

\begin{flushleft}

Joe Harris\\
Department of Mathematics\\
Harvard University\\
One Oxford Street\\
Cambridge, MA 02138\\
USA\\
{harris@math.harvard.edu}

\

Klaus Hulek\\
Institut f\"ur Mathematik\\
Universit\"at Hannover\\
Postfach 6009\\
D 30060 Hannover\\
Germany\\
{\tt hulek@math.uni-hannover.de}

\end{flushleft}

\end{document}